\documentclass[fleqn,12pt,twoside]{article}
\usepackage{amsmath}
\usepackage{amsfonts}

\usepackage{titletoc}
\usepackage{bbm}
\usepackage{txfonts}
\usepackage{espcrc1}

\title{\textbf{Existence of homoclinic solution in first order discrete Hamiltonian system} }
\author{Wenxiong Chen\address[]{School of Mathematical Science, Huaqiao University, Quanzhou 362021, P. R.
        China}
                      \thanks{
                          E-mail address:
                          chenwenxiong@amss.ac.cn. \hfil\break The author is supported by National Natural Science Foundation of China (No.11226115) and Natural Science Foundation of Huaqiao University (No.2014KJTD14)}
                           }
\begin{document}

\maketitle{\textbf{Abstract}}

\begin{abstract}

In this paper we consider the first order discrete Hamiltonian
system
$$\begin{cases}
x_1(n+1)-x_1(n)& =- H_{x_2}(n,x(n)), \\
                  x_2(n)-x_2(n-1)& =\ \  H_{x_1}(n,x(n)).
\end{cases}
$$ Where $n\in \mathbb{Z}$, $x(n)=$ $x_1 (n) \choose x_2 (n)$$ \in
\mathbb{R}^{2N}$, $ H(n,z)= \frac12S(n)z\cdot z +
 R(n,z)$ is periodic in $n$ and asymptotically quadratic as $|z| \to \infty$.  We will prove the existence
 of homoclonic solution by critical point theorem for strongly indefinite
functional.
\end{abstract}

{\bf Keywords:} discrete Hamiltonian system, homoclinic solution, asymptotically quadratic, variational method

\section{Introduction and main results}
In this paper we are interested in the following discrete first
order Hamiltonian system
$$
\begin{cases}
x_1(n+1)-x_1(n)& =- H_{x_2}(n,x(n)), \\
                  x_2(n)-x_2(n-1)& =H_{x_1}(n,x(n)).
\end{cases}\eqno(DHS)
$$
Where $n \in \mathbb{Z}$ and $x(n)=$ $x_1 (n) \choose x_2 (n)$ $\in \mathbb{R}^{2N}$,
$H(n,\cdot) \in \mathcal{C}^1 (\mathbb{R}^{2N}, \mathbb{R})$ depends
periodically on $n$ and has the form
$$
H(n,z)=\frac{1}{2} S(n)z \cdot z +R(n,z)
$$
 with $S(n)$ being a symmetric $2N \times
2N$ real matrix. Let $$
L x(n)=\binom{ x_1 (n+1) }{ x_2 (n)},\ \ \ \  \Delta x(n)=x(n+1)-x(n),$$ and $$\mathcal {J} = \left(\begin{array}{cc} 0&-I_N \\
I_N&0
\end{array} \right),$$ then we can rewrite system $(DHS)$ as follows
$$
\Delta Lx(n-1)= \mathcal {J} \nabla H(n,x(n)) \hskip0.5cm n \in
\mathbb{Z}. \eqno (DHS)^{\prime}
$$
We are interested in the existence of homoclinic solution
$x=(x(n))_{n\in \mathbb{Z}}$ of $(DHS)$, i.e.  $x \not\equiv 0$ and
$x(n) \to 0$ as $|n| \to \infty$.

System $(DHS)$ can be regarded as a discrete analog of continuous Hamiltonian system
$$
\begin{cases}
\dot x_1(t)& =- H_{x_2}(t,x(t)), \\
                 \dot x_2(t)& =H_{x_1}(t,x(t)).
\end{cases}\eqno(CHS)
$$
which have been largely studied in the literature
of the existence and multiplicity of homoclinic orbits by different approaches. Especially, there are some significant results for $(CHS)$ via variational method.
For details, we refer to [2-13] and references therein.

In last years, there have been many studies on discrete Hamiltonian systems by different approaches.
Abounding researches have been made on boundary value problems, oscillations
and asymptotic behavior of discrete Hamiltonian systems (see for example [20-24]). By critical point theory,
the existence and multiplicity of periodic solutions have been considered in [15-18]. It's well known Hamiltonian systems are very important in
the study of gas dynamics, fluid mechanics, relativistic mechanics and nuclear physics.
While it is well known that homoclinic solutions play an important role in analyzing the
chaos of Hamiltonian systems. If a system has the transversely intersected homoclinic
solutions, then it must be chaotic. If it has the smoothly connected homoclinic solutions, then it cannot stand the perturbation, its perturbed system probably produces chaotic
phenomena. Therefore, it is of practical importance and mathematical significance to
consider the existence of homoclinic solutions of Hamiltonian systems. As we know, there are not much results for homoclinics of discrete Hamiltonian systems.
 As in \cite{D & L}, the existence of homoclinic solutions has been obtained in second
order discrete Hamiltonian systems. Recently, Chen, Yang and Ding \cite{C&Y&D} proves existence and multiplicity of homoclinics in first order Hamiltonian systems with the nonlinear term being super-quadratic. Hence our aim of this paper is to establish some existence results for first order discrete Hamiltonian system with asymptotically quadratic term.

We will discuss our results variationally. Inspired by [11], we discuss the associated linear self-adjoint operator $A+S$(defined in Section 2). By the spectrum of $A+S$, we establish he variational framework for $(DHS)$. And we show that the  associated functional $\Phi$ is strongly indefinite.

To state our results, we use the notation
$\mathcal{J}_0=\left(\begin{array}{cc} 0&-I_N \\-I_N&0\end{array} \right)$. For a sequence of symmetric matrixs $\{L(n)\}$,
let $\mathfrak{e}(n)$ be the set of all eigenvalues of $L(n)$ and set $$\lambda_L =\inf_n \mathfrak{e}(n), \hskip0.5cm \Lambda_L =\sup_n \mathfrak{e}(n).$$
In particular, we set $\lambda_0:= \lambda_{\mathcal{J}_0 S}$ and $\Lambda_0:= \Lambda_{\mathcal{J}_0 S}$. And we will use the notation
$$\tilde R(n,z)=\frac{1}{2} \nabla R(n,z)z-R(n,z)$$
 We make the following
hypotheses:
\begin{itemize}
\item[$(R_0)$] There is a positive integer $T$ such that $S(n+T)=S(n)$, and
$\mathcal{J}_0 S(n)$ is symmetric and positive definite for all $n\in\mathbb{Z}$.
\item[$(R_1)$] $R(n+T,z)=R(n,z)$, $\forall n\in\mathbb{Z}$ $\forall
z\in\mathbb{R}^{2N}$, $R(n,\cdot) \in
\mathcal{C}^1(\mathbb{R}^{2N},\mathbb{R})$.
\item[$(R_2)$] $R(n,z)\geq 0$, $\nabla R(n,z)=o(|z|)$ as $|z| \to 0$.
\item[$(R_3)$] $\nabla R(n,z)-S_{\infty}(n)z=o(|z|)$ as
 $|z|\to\infty$, where $S_{\infty}(n)$ is a symmetric matrix with
 $\lambda_{\infty}:=\lambda_{S_{\infty}} >2+\Lambda_0$
\item[$(R_4)$]$\tilde R(n,z)\geq0$ and there
is $\delta_0\in (0,\lambda_0)$ such that if $|\nabla
R(n,z)|\geq(\lambda_0-\delta_0)|z|$ then $\tilde R(n,z)\geq \delta_0$.
\end{itemize}

Then we have the main result of this paper.

 {\bf Theorem 1.1} Let
$(R_0)-(R_4)$ be satisfied. Then the system $(DHS)$ has at least one homoclinic
solution.

This paper is organized as follows. In Section 2, we establish the variational framework of the problem
and recall some abstract critical point theories on strongly indefinite functional.
In Section 3, we discuss the linking structure of $\Phi$ and the behavior of $(C)_c$-sequence. Finally, in Section 4, we  prove our results.

\section{Variational setting}
Let $E:=l^2(\mathbb{Z},\mathbb{R}^{2N})$. $E$ is a Hilbert space
with the usual inner product and norm
$$(x,y)_{l^2}=\sum_n x(n)
\cdot y(n) \hskip0.5cm |x|_{l^2}^2=\sum_n |x(n)|^2 \hskip0.5cm x,y
\in E \eqno(2.1)$$
On $E$
we define functional $\Phi$, for any $x\in E$,
$$\Phi(x)=-\frac{1}{2}\sum_n \mathcal {J}\Delta L x(n-1) \cdot x(n) -
\frac{1}{2} \sum_n S(n)x(n) \cdot x(n)- \sum_n R(n,x(n)). \eqno
(2.2)$$

For convenience, we define operators as follows.
$$A: E \to E: \hskip0.5cm Ax=(z(n))_{n\in \mathbb {Z}} \hskip0.5cm z(n)=-\mathcal
{J}\Delta Lx(n-1),\hskip0.5cm x\in E \eqno(2.3)$$ $$S: E \to E: \hskip0.5cm
\hskip0.5cm Sx=(z(n))_{n\in \mathbb {Z}} \hskip0.5cm z(n)=-S(n)x(n),
\hskip0.5cm x\in E.\eqno(2.4)$$Then $A$ and $S$
are linear bounded self-adjoint operators(see\cite{C&Y&D} ).

Moveover, we set $$\Psi(x):=\sum_n R(n,x(n))\eqno (2.5)$$
Thus, we can rewrite functional $\Phi$:
$$\Phi(x)=\frac{1}{2}
((A+S)x,x)_{l^2}-\Psi(x). \eqno(2.6)$$

Since Hilbert space $E=l^2$ embeds continuously into $l^p
(\mathbb{Z},\mathbb{R}^{2N})$($2<p\leq \infty$), i.e.
$$|x|_{l^p} \leq |x|_{l^2} \hskip1cm x\in E,\eqno(2.7)$$ and $(R_1)-(R_3)$ imply that
, for any $\varepsilon>0$,$p\geqslant 2$, there is $C_{\varepsilon} > 0$ such that
$$|\nabla R(n,z)| \leq \varepsilon |z|+ C_{\varepsilon} |z|^{p-1} \eqno (2.8)$$ and
$$|R(n,z)| \leq \varepsilon |z|^2+ C_{\varepsilon} |z|^p. \eqno
(2.9)$$
Hence the functional $\Phi$ is well defined. By a standard argument, one can obtain that $\Phi$ is $\mathcal{C}^1$ and has
 Fr\'{e}chet derivative of the
following form
$$\Phi^{\prime}(x)y= -\sum_n \mathcal {J}\Delta Lx(n-1) \cdot y(n) -
\sum_n S(n)x(n) \cdot y(n)- \sum_n \nabla R(n,x(n))\cdot y(n) \eqno
(2.10)$$
for $x,y \in E$.

Obviously, if $x \in E$ is a critical point of $\Phi$ then $x$
is a solution of (1), moreover, $x(n) \to 0$ as $|n| \to \infty$. Since $S(n)$ and $R(n,z)$ are $T$-periodic on $n$,
 it is not difficult to see that $\Phi$ is $T$-periodic.
\\

In order to establish a variantional setting for the system $(DHS)$
 we study the spectrum of the associated linear self-adjoint operator $A+S$.
Let $\sigma (A+S)$, $\sigma_e (A+S)$ denote, respectively, the spectrum and essential spectrum of $A+S$.
Supposing $(R_0)$
holds, by the definition of $\lambda_0$, it is easy to see that $\lambda_0>0$.
\\

{\bf Proposition 2.1} Assume $(R_0)$ is satisfied. Then
\begin{itemize}
\item[$1^{\textordmasculine}$] $\sigma (A+S)=\sigma_e (A+S)$;
\item[$2^{\textordmasculine}$] $\sigma (A+S) \cap (0,\infty)\neq \O$
and $\sigma (A+S) \cap (-\infty,0)\neq \O$;
\item[$3^{\textordmasculine}$] $\sigma (A+S) \subset
[-\Lambda_0- 2,-\lambda_0] \bigcup [\lambda_0,\Lambda_0+2]$.
\end{itemize}

Proof. \hskip0.5cm The proofs of $1^{\textordmasculine}$,$2^{\textordmasculine}$ and $\sigma (A+S) \subset \mathbb{R}\setminus
(-\lambda_0,\lambda_0)$, one can obtain in \cite{C&Y&D}. To prove that $3^{\textordmasculine}$ holds, it is sufficient
to show that $\|A+S\|\leq 2 +\Lambda_0$. Note that
\begin{eqnarray*}
|Ax|^{2}_{l^2} &=&(Ax,Ax)_{l^2}=\sum_n  (-\mathcal
{J}\Delta Lx(n-1)) \cdot (-\mathcal
{J}\Delta Lx(n-1)) \\&=& \sum_n \{|x_2 (n)|^2 +|x_2 (n-1)|^2 +|x_1 (n)|^2+|x_1 (n+1)|^2  \\& & \hskip0.6cm -2x_2 (n) \cdot x_2 (n-1)-2x_1 (n) \cdot x_1 (n+1)\}
\\ &\leq& \sum_n \{2|x_1 (n)|^2 +2|x_1 (n+1)|^2 +2|x_2 (n)|^2 +2|x_2 (n-1)|^2\}
\\&=& 4|x|^2 _{l^2}
\end{eqnarray*}
Observe that $\mathcal{J}^2_0=I$ and $\mathcal{J}_0 S(n)=S(n)\mathcal{J}_0$, we have
\begin{eqnarray*}
|Sx|^2_{l^2}&=& (Sx,Sx)_{l^2}=\sum_n (-S(n)x(n)) \cdot (-S(n)x(n)) =\sum_n \mathcal{J}^2_0 S(n)x(n) \cdot \mathcal{J}^2_0 S(n)x(n)
\\&=& \sum_n \mathcal{J}_0 S(n)\mathcal{J}_0x(n) \cdot \mathcal{J}_0 S(n)\mathcal{J}_0x(n)
=\sum_n (\mathcal{J}_0 S(n))^2 \mathcal{J}_0x(n) \cdot \mathcal{J}_0x(n).
\end{eqnarray*}By definition of $\Lambda_0$, we have that $|Sx|^2_{l^2}\leq \Lambda^2_0 |x|^2_{l^2}$. It follows that $\|A\|\leq 2$ and $\|S\|\leq \Lambda_0$.
Therefore, $\|A+S\|\leq \|A\|+\|S\| \leq 2+\Lambda_0$. The proof is complete.

Due to the spectrum of $A+S$, $E=l^2(\mathbb{Z},\mathbb{R}^{2N})$ possesses the orthogonal decomposition $$E=E^- \oplus E^+ ,
\hskip0.5cm x= x^- + x^+, \eqno (2.11)$$ corresponding to the
spectrum decomposition of $A+S$ such that $$((A+S)x,x)_{l^2}\leq
-\lambda_0 |x|_{l^2}^2 \hskip5pt on \hskip5pt E^- \hskip5pt and
\hskip5pt ((A+S)x,x)_{l^2}\geq\lambda_0 |x|_{l^2}^2 \hskip5pt on
\hskip5pt E^+ .\eqno (2.12)$$
Let $|A+S|$ denote the absolute value of $A+S$, we equip $E$ with the
inner product
$$( x,y ) =(|A+S|^{1/2}x,|A+S|^{1/2}y)_{l^2}.$$ Then $(E,(\cdot,\cdot))$ is a Hilbert space, and it has the associated norm
 $\|x\|= (x, x)^{1/2}$. By $3^{\textordmasculine}$ of Proposition 2.1, one can obtain that
$$\lambda_0 |x|_{l^2}\leq \|x\|\leq (2+\Lambda_0)|x|_{l^2}, \eqno(2.13)$$
which implies that $(E,\|\cdot\|)$ is equivalent to
$(E,|\cdot|_{l^2})$. It is not difficult to see that the
decomposition of $E$ is orthogonal with respect to both $(\cdot,
\cdot)_{l^2}$ and $( \cdot,\cdot )$. Now, we can rewrite the functional $\Phi$ as
$$\Phi(x)=\frac{1}{2}\|x^+\|^2-\frac{1}{2}\|x^-\|^2-\Psi(x).\eqno(2.14)$$
Hence, by the Proposition 2.1, $\Phi$ is strongly indefinite.

To study the critical point of $\Phi$, we recall some abstract critical point theory developed
in \cite{T & D1}.

Let $Z$ be a Banach space with direct sum decomposition $Z=X\oplus
Y$ and corresponding projections $P_X$, $P_Y$ onto $X$, $Y$,
 respectively. For a functional $\Phi \in \mathcal{C}^1(Z,\mathbb{R})$ we write
$\Phi_a = \{z\in Z : \Phi (z)\geq a\}$, $\Phi^b = \{z\in Z : \Phi (z) \leq b\}$
and $\Phi^b_a =\Phi_a \cap  \Phi^b$. Recall that $\Phi$ is said to be weakly sequentially lower semicontinuous
if for any $z_n \rightharpoonup z$ in $Z$ one has $\Phi(z) \leq \liminf_{n \to \infty} \Phi(z_n)$, and $\Phi^{\prime}$ is
said to be weakly sequentially continuous if $\lim_{n \to \infty} \Phi^{\prime}(z_n)w = \Phi^{\prime}(z)w$ for each $w \in Z$.
A sequence $(z_n) \in Z$ is
said to be a $(C)_c$-sequence if $\Phi (z_n) \to c$ and $\Phi^\prime
(z_n)\to 0$. $\Phi$ is said to satisfy the $(C)_c$-condition if any $(C)_c$-sequence has a convergent subsequence.

From now on, let $X$ be separable and reflexive, and fix a
countable dense subset $\mathcal{S} \subset X^*$. For each $s\in \mathcal{S}$ there is a semi-norm on $Z$
defined by $$p_s : Z \to \mathbb{R}, p_s(z) =
|s(x)| + \|y \| \hskip0.5cm for \hskip0.3cm z = x + y \in Z=X\oplus
Y.$$ We denote by $\mathcal{T}_\mathcal{S}$ the induced topology,
Let $w^*$ denote the $weak^*$-topology on $Z^*$.

Suppose:
\begin{itemize}
\item[$(\Phi_0)$] For any $c \in \mathbb{R}$,
 $\Phi_c$ is $\mathcal{T}_\mathcal{S}$-closed, and $\Phi
^\prime:(\Phi_c,\mathcal{T}_\mathcal{S}) \to (E^*,w^*)$ is
 continuous.
\item[$(\Phi_1)$] For any $c>0$, there exists $\zeta>0$ such that $\|z\| \leq \zeta \|P_Y z\|$ for all $z\in \Phi_c$.
\item[$(\Phi_2)$] There exists $\rho >0$ with $\kappa =\inf \Phi(S_\rho Y)>0$ where $S_\rho Y=\{ y \in Y : \|y \|=\rho
\}$.
\end{itemize}
Then the following theorem is a special case of Theorem 4.4 of \cite{T & D1}.

{\bf Theorem 2.2} Let $(\Phi_0)-(\Phi_2)$ be satisfied and assume
there is $R>\rho$ and $e\in Y$ with $\|e\|=1$ such that $\sup
\Phi(\partial Q ) \leq \kappa$ where $Q=\{z=x+te: t\geq 0, x\in X, \|z\|<R
\}$. Then there is a $(C)_c$-sequence for $\Phi$ with $c\in
[\kappa, \sup \Phi(Q)]$.

To check that functional $\Phi$ satisfies $(\Phi_0)$, the following proposition (cf. [1]) is a key tool. 

{\bf Proposition 2.3}  Suppose $\Phi \in C^1(Z,\mathbb{R})$ is of the
form $$\Phi (z) = \frac{1} {2}\|y\|^2-\frac{1}{2}\|x\|^2-\Psi (z)
\hskip0.5cm z = x + y \in Z= X \oplus Y \eqno (2.15)$$ such that
\begin{itemize}
 \item[(i)] $\Psi \in C^1(Z,\mathbb{R})$ is bounded from below.
 \item[(ii)] $\Psi$ is weakly
sequentially lower semicontinuous
\item[(iii)]   $\Psi^\prime $ is weakly sequentially
continuous.
\item[(iv)] $\nu : Z \to \mathbb{R}, \nu (z) = \| z\|^2$,
is $C^1$ and $\nu^\prime: (Z, \mathcal T_w) \to (Z^*, \mathcal
T_{w^*} )$ is sequentially continuous.
\end{itemize}
Then $\Phi$ in (2.15) satisfies $(\Phi_0)$.

\section{Linking structure and $(C)_c$ sequence}

In order to apply Theorem 2.1, we study the linking structure of $\Phi$.

{\bf Lemma 3.1}  Let $(R_0)-(R_3)$ be satisfied. Then there exists $\rho >0$ such that
$\kappa =\inf \Phi(S_\rho ^+)>0$ where $S_\rho ^+=\partial B_\rho \bigcap E^+$.

Proof \hskip0.5cm Choose $p>2$ such that (2.9) holds for any
$\varepsilon>0$. We get that
$$\Psi (x)\leq \varepsilon |x|^2_{l^2}+C_\varepsilon |x|^p_{l^p}\leq (2+\Lambda_0)(\varepsilon \|x\|^2+C_\varepsilon \|x\|^p)$$
for all $x \in E$. We choose $\varepsilon$ small enough, then the lemma holds from the form of $\Phi$.
\\

{\bf Lemma 3.2}  Let $(R_0)-(R_3)$ be satisfied.
Then for any $e\in E^+$ with $\|e\|=1$
 there exists $R>0$ such that $\Phi
(x)\leq \kappa$ for all $x \in E^-\oplus \mathbb{R}e$, with $\|x\| \geq R$.

Proof \hskip0.5cm  It is sufficient to
show that $\Phi (x) \to -\infty$ as $x \in  E^-\oplus \mathbb{R}e$, $\|x\| \to \infty$.
 Arguing by indirectly, assume that for some sequence $x_k=s_k e+x^-_k \in E^-\oplus \mathbb{R}e$
with $\|x_k\| \to \infty$, there is $M>0$ such that $\Phi (x_k)\geq
-M$ for all $k$. Then, setting $y_k= x_k/\|x_k\|:=t_k e+ y^-_k $, we have
$\|y_k\|=1$, $y_k\rightharpoonup y$, $y^- _k \rightharpoonup y^-$,
$t_k \to t \in \mathbb{R}$ and
$$-\frac{M} {\|x_k\|^2}\leq \frac{1} {2}|t_k|^2-\frac{1} {2}\|y^-_k\|^2
- \frac{\sum_n R(n,x_k(n))} {\|x_k\|^2} . \eqno (3.1)$$ Remark that
$t\neq 0$. Indeed, if not then it follows from (3.1) that
$$0\leq \frac{1} {2}\|y^-_k\|^2
+\frac{\sum_n R(n,x_k(n))} {\|x_k\|^2} \leq \frac{1} {2}|t_k|^2
+\frac{M} {\|x_k\|^2},\eqno(3.2)$$ in particular,
 $\|y^-_k\| \to 0$, hence $1=\|y_k\| \to 0$, a contradiction.

Since $(R_3)$, there holds
\begin{eqnarray*}
|t|^2 - \|y^-\|^2 -\sum_n S_{\infty}(n)
y(n)\cdot y(n)&\leq&
 \|te\|^2 - \|y^-\|^2 -\lambda_{\infty} |y|^2_{l^2}
  \\&\leq&(2+\Lambda_0)|t|^2|e|^2_{l^2} - \|y^-\|^2 -\lambda_{\infty} |t|^2|e|^2_{l^2}-\lambda_{\infty} |y^-|^2_{l^2}
 \\&=&(2+\Lambda_0-\lambda_{\infty})|t|^2|e|^2_{l^2}- \|y^-\|^2-\lambda_{\infty} |y^-|^2_{l^2}<0.
 \end{eqnarray*}
 Hence
  for some $\tilde N>0$
  $$|t|^2 - \|y^-\|^2 -\sum _{|n|\leq\tilde N}  S_{\infty}(n)
y(n)\cdot y(n) <0 .\eqno (3.3)$$
  Let $G(n,z):=R(n,z)-\frac{1}{2}S_{\infty}(n)z\cdot z$, by $(R_3)$ and (2.9), one gets
  that $|G(n,z)| \leq C_2 |z|^2$.

  Since $$\lim_{k \to \infty} \Big ( \sum _{|n|\leq\tilde N} \frac{R(n,x_k(n))} {\|x_k\|^2}
  -\frac{1} {2} \sum _{|n|\leq\tilde N}  S_{\infty}(n)
y_k(n)\cdot y_k(n)\Big) =
  \lim_{k \to \infty} \sum _{|n|\leq\tilde N} \frac{G(n,x_k(n))} {\|x_k\|^2}
  $$$$=
  \lim_{k \to \infty} \sum _{|n|\leq\tilde N} \frac{ G (n,x_k(n)) |y_k(n)|^2} {|x_k(n)|^2}, $$
 and we have that $$\Big| \sum _{|n|\leq\tilde N} \frac{ G (n,x_k(n)) |y_k(n)|^2} {|x_k(n)|^2} \Big|
   \leq  \sum _{|n|\leq\tilde N} \frac{ |G (n,x_k(n))| |y_k(n)|^2}
{|x_k(n)|^2}.$$ For $ |n| \leq \tilde N$, if $y_k(n) \to 0$ then $
\frac{ |G (n,x_k(n))| |y_k(n)|^2} {|x_k(n)|^2}\leq C_2 |y_k(n)|^2
\to 0$, otherwise, if $y_k(n) \not \to 0$, then $|x_k(n)|\to
\infty$, which yields that $\frac{ |G
(n,x_k(n))| |y_k(n)|^2} {|x_k(n)|^2} \leq C_3 \frac{ |G (n,x_k(n))|
} {|x_k(n)|^2}\to 0 $ $(C_3>0)$. Thus,
$$\sum _{|n|\leq\tilde N} \frac{ G (n,x_k(n)) |y_k(n)|^2} {|x_k(n)|^2} \to 0 .$$

  It follows that $$0\leq \lim_{k \to \infty} \Big ( \frac{1} {2} |t_k |^2-\frac{1} {2}\|y^-_k\|^2
- \sum _{|n|\leq\tilde N} \frac{R(n,x_k(n))} {\|x_k\|^2}\Big )
   \leq \frac{1} {2} \Big( |t|^2 - \|y^-\|^2 -\sum _{|n|\leq\tilde N}  S_{\infty}(n)
y(n)\cdot y(n) \Big) <0, $$
that is a contradiction.

It follows from Lemma 3.1 and Lemma 3.2 that $\Phi$ has linking structure which is showed by following lemma.

{\bf Lemma 3.3} Under the assumptions of Lemma 3.2, letting $e\in E^+$ with $\|e\|=1$, there is $R_0>\rho$ such that $\sup(\partial Q)\leq\kappa$ 
where $Q:=\{x=x^-+te: t\geq 0, x^-\in E^-, \|x\|<R_0
\}$.

\

\
Now we discuss the behavior of $(C)_c$-sequence.

{\bf Lemma 3.4} Let $(R_0)-(R_4)$ be satisfied, then any $(C)_c$-sequence of $\Phi$ is bounded.

Proof \hskip0.5cm Let $(x_k) \subset E$ be such that
 $$\Phi (x_k)  \to c \hskip0.5cm  \hskip0.5cm (1+\|x_k\|)\Phi ^\prime (x_k) \to 0.  \eqno (3.4)$$
Then, for some $C_0>0$,
 $$C_0\geq \Phi (x_k)- \frac{1} {2} \Phi ^\prime (x_k) x_k=
\sum_n \tilde R(n,x_k(n)). \eqno (3.5)$$ Arguing by indirectly
assume up to a subsequence $\|x_k \| \to \infty$. Set $y_k = x_k
/\|x_k\| $. Then, $\|y_k\|=1$. Remark that
$$\Phi ^\prime (x_k)(x^+ _k- x^- _k) =
\|x_k\|^2 \Big ( 1- \frac{\sum_n \nabla R(n,u_k(n)(y^+_k(n)-y^-
_k(n))} {\|x_k\|} \Big)$$ it follows from (3.4) that
$$\frac{\sum_n \nabla R(n,u_k(n)(y^+_k(n)-y^-
_k(n))} {\|x_k\|} \to 1 \eqno (3.6)$$

If $ \exists \tau >0$, $ \exists n_k \in \mathbb{Z}$ such that
$|y_k(n_k)|\geq \tau$. There exist $l_k \in \mathbb{Z}$ such that $
0 \leq n_k -l_k T \leq T-1$. Set $\tilde x _k =\{ \tilde x_k(n) \}$
where $\tilde x_k(n)= x_k(n+l_k T)$ and $\tilde y _k =\{ \tilde
y_k(n)\}$ where $\tilde y_k(n)= y_k(n+l_k T)$. For any $w \in E$
setting $\tilde w _k =\{ \tilde w_k(n) \}$ where $\tilde w_k(n)=
w_k(n -l_kT)$, and define operator $S_{\infty}$ as

$$S_{\infty}: E \to E: \hskip0.5cm
\hskip0.5cm S_{\infty}x=(z(n))_{n\in \mathbb {Z}} \hskip0.5cm z(n)=S_{\infty}(n)x(n),
\hskip0.5cm x\in E.$$
 We have that
$$\Phi^ \prime (x_k) \tilde w_k=(x^+ _k -x^- _k,\tilde w_k )-(S_{\infty} x_k,\tilde w_k)_{l^2}- \sum_n \nabla R_(n,x_k(n))\tilde w_k(n)$$
 $$=\|x_k\|\Big((y^+ _k -y^- _k,\tilde w_k ) -(S_{\infty} y_k, \tilde w_k)_{l^2}-
 \sum_n \nabla  R(n,x_k(n)) \tilde w_k(n) \frac{|y_k(n)|}{|x_k(n)|} \Big )$$
By $A+S$, $S_{\infty}$ and $\nabla R(n,z)$ are periodic in $n$, then
 $$\Phi^ \prime (x_k) w_k=\|x_k\|\Big(( \tilde y^+ _k -\tilde y^- _k,w ) -(S_{\infty} \tilde y_k,w)_{l^2}-
 \sum_n  \nabla  R(n,\tilde x_k(n)) w(n) \frac{|\tilde y_k(n)|}{|\tilde x_k(n)|} \Big )$$
 This follows
 $$(\tilde y^+ _k -\tilde y^- _k,w ) -(S_{\infty} \tilde y_k,w)_{l^2}-
 \sum_n  \nabla  R(n,\tilde x_k(n)) w(n) \frac{|\tilde y_k(n)|}{|\tilde x_k(n)|}\to 0$$
 Since $\|\tilde y_k\|=\|y_k\|=1$, then  up to subsequence we obtain that $\tilde y_k
 \rightharpoonup \tilde y$, $\tilde y_k(n) \to \tilde y_(n)$. For some $n_0: 0 \leq n_0 \leq (T-1)$ such that $|\tilde
 y_k(n_0)|\geq \tau$, hence, $\tilde y \neq 0$. By (2.8), one
 gets that $|\nabla R(n,z)| \leq C_4 |z|$. For $w \in E=l^2$ and $\varepsilon >0$,
 there exists $\tilde N \in
 \mathbb{N}$ such that $\sum_{|n|> \tilde N} |w(n)|^2 <
\varepsilon^2$. Note that
 $$\Big| \sum_n  \nabla  R(n,\tilde x_k(n)) w(n) \frac{|\tilde y_k(n)|}{|\tilde x_k(n)|}
 \Big| \leq \sum_n
 | R(n,\tilde x_k(n))| |w(n)| \frac{|\tilde y_k(n)|}{|\tilde x_k(n)|}
$$
$$ \leq \sum_{|n|\leq \tilde N} |\nabla R(n,\tilde x_k(n))| |w(n)| \frac{|\tilde y_k(n)|}{|\tilde x_k(n)|} +
\sum_{|n|>\tilde N} |\nabla R(n,\tilde x_k(n))| |w(n)| \frac{|\tilde y_k(n)|}{|\tilde x_k(n)|}$$

$$\leq \sum_{|n|\leq \tilde N} |\nabla R(n,\tilde x_k(n))| |w(n)| \frac{|\tilde y_k(n)|}{|\tilde x_k(n)|} +
C_4 |y_k|_{l^2} (\sum_{|n|>\tilde N} |w(n)|^2)^{\frac{1}{2}} $$
$$ \leq
\sum_{|n|\leq \tilde N} |\nabla R(n,\tilde x_k(n))| |w(n)| \frac{|\tilde y_k(n)|}{|\tilde x_k(n)|} + C_5 \varepsilon.$$
Similarly to the proof of Lemma 3.2, we can obtain that
$$\sum_{|n|\leq \tilde N} |\nabla R(n,\tilde x_k(n))| |w(n)| \frac{|\tilde y_k(n)|}{|\tilde x_k(n)|} \to 0.$$
  It follows that
  $$\sum_n |\nabla R(n,\tilde x_k(n))| |w(n)| \frac{|\tilde y_k(n)|}{|\tilde x_k(n)|} \to
  0,$$
  hence $$(\tilde y^+ -\tilde y^-, w ) -(S_{\infty} \tilde y , w)_{l^2}=0.$$
  Thus, $0$ is an eigenvalue of the operator $A+S-S_{\infty}$ with eigenfunction $\tilde y$. We claim that is impossible.
  Indeed, by the definition of $\lambda_{\infty}$ and (2.13), one has that, for any $x \in E$
  $$|(A+S)x-S_{\infty}x|_{l^2}\geq |S_{\infty}x|_{l^2}-|(A+S)x|_{l^2}\geq (\lambda_{\infty}-2-\Lambda_0)|x|_{l^2}.$$
It follows from $(R_3)$ that $\lambda_{\infty}-2-\Lambda_0>0$ which yields that $0\notin \sigma (A+S-S_{\infty})$.

  If for any $\tau >0$ and $ n$ there holds $|y_k(n)|<\tau$, then
  $|y_k|_{l^{\infty}} \to 0$. Since  $|y_k|^p_{l^p} \leq |y_k|^{p-2}_{l^{\infty}}
  |y_k|^2_{l^2}$ $(p>2)$
  and $y_k$ is bounded in $E=l^2$, we get that $|y_k|_{l^p} \to 0 $ $(p>2)$.
  In virtue of $(R_4)$,  we set
  $$I_k := \left\{  n \in \mathbb{Z} : \frac{|\nabla R(n, x_k(n))|} {|x_k(n)|}  \leq \lambda_0 - \delta_0  \right\}$$
  Since $\lambda_0 |y_k|^2_{l^2} \leq \|y_k\|^2=1$, we have
  \begin{eqnarray*}
  \Big| \sum _{I_k} \frac{ \nabla R(n, x_k(n))(y^+_k(n)-y^- _k(n))}
{\|x_k\|} \Big| &=& \Big| \sum _{I_k} \frac{
\nabla R(n, x_k(n))(y^+_k(n)-y^-_k(n))|y_k(n)|}{|x_k(n)|} \Big|
\\&\leq&
(\lambda_0-\delta_0) |y_k|^2_{l^2} \leq \frac{\lambda_0-\delta_0}{\lambda_0} < 1
\end{eqnarray*}
 for all $k$. Let
$I^c_k =\mathbb{Z} \backslash I_k$, jointly with (3.6), implies
that
$$\lim_{k \to \infty} \Big| \sum _{I^c_k} \frac{ \nabla R(n, x_k(n))(y^+_k(n)-y^- _k(n))}
{\|x_k\|} \Big| >1- \frac{\lambda_0-\delta_0}{\lambda_0} = \frac{\delta_0}
{\lambda_0}.$$
 Recalling (2.8), we can choose $C>0$ such that $|\nabla R(n ,z)|\leq
C|z|$, there holds for an arbitrarily fixed $s>2$,
\begin{eqnarray*}\Big| \sum _{I^c_k} \frac{ \nabla R(n, x_k(n))(y^+_k(n)-y^- _k(n))}
{\|x_k\|} \Big| &\leq& C \sum _{I^c_k} |y^+ _k(n) -y^- _k(n)|
|y_k(n)| \leq C|y_k|_{l^2}|I^c_k|^{(s-2)/s} |y_k|_{l^s}
\\&\leq& \frac
{C}{\lambda_0} |I^c_k|^{(s-2)/s} |y_k|_{l^s}\
\end{eqnarray*}
Since $|y_k|_{l^s} \to 0$, one gets that $|I^c_k| \to \infty$. By
 $(R_4)$, $\tilde R(n, x_k(n))\geq \delta_0$ on $I^c_k$, hence
 $$\sum_n \tilde R(n, x_k(n)) \geq \sum_{I^c_k} \tilde R(n, x_k(n)) \geq |I^c_k| \delta_0 \to \infty$$
 contrary to (3.5). The proof is finished.

\section{Proof of main result}
We are now in a position to give the proof of our main result. In order to apply the abstract Theorem 2.2,  we
choose $X=E^-$ and $Y=E^+$ with $E^\pm$ given in Section 2. $X$ is separable and reflexive and let $\mathcal{S}$
be a countable dense subset of $X^*$. First we have

{\bf Lemma 4.1}  $\Phi$ satisfies $(\Phi_0)$ and $(\Phi_1)$.

Proof. \hskip0.5cm In virtue of the form of $\Phi$ and Proposition 2.3, to show that $\Phi$ satisfies $(\Phi_0)$ it is sufficient to show that $\Psi$ is bounded from below,
$\Psi$ is weakly sequentially lower semicontinuous and $\Psi^{\prime}$ is weakly sequentially continuous.

 Firstly, since $R(n,z)$ is non-negative, so is $\Psi$. Secondly, let $x_k
\rightharpoonup x$ in $E$. We have that $x_k (n) \to x(n)$ as $k \to
\infty$, hence, $R(n,x_k(n)) \to R(n,x(n))$.
 Thus,
$$\Psi (x) =\sum_n \lim_{k\to \infty}R(n,x_k(n))
\leq \liminf_{k\to \infty} \sum_n R(n,x_k(n)) = \liminf_{k\to
\infty} \Phi (x_k)$$
which implies that $\Psi$ is weakly sequentially lower semicontinuous.

Thirdly, let
$x_k \rightharpoonup x$  in $E$. By (2.8), we choose $C_1>0$ such that
$|\nabla R(n,z)| \leq C_1 |z|$. For any $y \in E$, one can get that
for any $\varepsilon>0$ there is $N\in \mathbb{N}$ such that
$\sum_{|n|\geq M} |y(n)|^2 <\varepsilon$. By H\"{o}lder inequality,
\begin{eqnarray*}|\Psi^{\prime}(x_k)y-\Psi^{\prime}(x)y|&\leq& \Big|\sum_{|n|> N}
(\nabla R(n,x_k(n)-\nabla R(n,x(n))\cdot y(n) \Big| \\ &
&+\Big|\sum_{|n|\leq N} (\nabla R(n,x_k(n)-\nabla R(n,x(n))\cdot
y(n) \Big| \\ &\leq& C_1 \varepsilon (|x_k|_{l^2}+|x|_{l^2}
)+\Big|\sum_{|n|\leq N} (\nabla R(n,x_k(n)-\nabla R(n,x(n))\cdot
y(n) \Big| .\end{eqnarray*}  Observe that $\nabla R(n,x_k(n))\to
\nabla R(n,x(n))$ for each $n\in \mathbb{Z}$. It follows from
$|n|\leq N$ is finite that $\Big|\sum_{|n|\leq N} (\nabla
R(n,x_k(n)-\nabla R(n,x(n))\cdot y(n) \Big| \to 0$. Therefore we
obtain that $\Psi^{\prime}$ is weakly sequentially continuous.

For any $c>0$,since $\Psi$ is non-negative, it is not difficult to see that $\|x^-\|\leq \|x^+\|$.
 Hence one can get that $\|x\|\leq 2\|x^+\|$ for $x\in \Phi_c$ which shows that $(\Phi_1)$ holds. The proof is complete.
 \\
{\bf Proof of therem 1.1} \hskip0.5cm
 Lemma 4.1 implies that $\Phi$ satisfies $(\Phi _0)$ and  $(\Phi
_1)$. And Lemma 3.3 shows that $\Phi$ has linking
structure. Hence there is a $(C)_c$-sequence $(x_k)$ with level
$c\geq \kappa >0$. Lemma 3.4 shows that $(x_k)$ is bounded: $\|x_k
\|\leq M$. In addition, $$c=\lim_{k \to \infty} \Big (\Phi (x_k)-
\frac{1} {2} \Phi ^\prime (x_k) x_k \Big )= \lim_{k \to \infty} \sum_n
\tilde R(n,x_k(n)) \eqno (4.1)$$
Remark that there are $\tau >0$ and $n_k \in \mathbb{Z}$ such that
$|x_k(n_k)| \geq \tau$. Indeed, if not, it follows from the proof
of Lemma 3.3 that $|x_k|_{l^p} \to 0 (p>2)$. By (2.8) and (2.9), choose
$p>2$, such that for any $\varepsilon >0$, there is $C_{\varepsilon}
>0$ satisfying $\tilde R(n ,z) \leq \varepsilon \lambda^2_0 M^{-2} |z|^2
+ C_{\varepsilon} |z|^p$, then it follows from (4.1) that, for
$\varepsilon <c$,
$$c=\lim_{k \to \infty} \sum_n \tilde R(n,x_k(n)) \leq  \lim_{k \to \infty} \Big (
   \varepsilon\lambda^2_0  M^{-2} |x_k|^2_{l^2} + |x_k|^p_{l^p}  \Big) \leq
   \varepsilon,$$ a contradiction.
 Make proper shifts similar
to the proof of Lemma 3.4, then passing to a subsequence $\tilde
x_{k}$ such that there is $\tau >0$ and $ 0 \leq n_0 \leq (T-1)$
independent of $k$,  $|\tilde x_k(n_0)| \geq \tau$. Due to periodicity of the coefficients, $\tilde x_k$ is also a
Cerami-sequence for $\Phi$ at c. Hence $\tilde
x_k \rightharpoonup \tilde x \neq 0$, thus, we obtain a nontrivial
critical point $\tilde x$ of $\Phi$.

\end{document}